\documentclass[a4paper,10pt]{article}
\usepackage{amssymb,amsmath}

\newcommand{\Z}{\mathbb{Z}}
\newcommand{\R}{\mathbb{R}}

\newcommand{\twosum}[2]{\sum_{\substack{#1\\#2}}}
\newcommand{\ep}{\varepsilon}
\newcommand{\lk}{\lambda_k}
\newcommand{\bal}{\boldsymbol{\alpha}}
\newcommand{\beql}[1]{\begin{equation}\label{#1}}
\newcommand{\eeq}{\end{equation}}

\newtheorem{theorem}{Theorem}
\newtheorem{corollary}{Corollary}
\newtheorem{lemma}{Lemma}

\begin{document}
\title{A New $k$-th Derivative Estimate for Exponential Sums via
  Vinogradov's Mean Value}
\author{D.R. Heath-Brown\\Mathematical Institute, Oxford}
\date{}
\maketitle

\begin{flushright}
{\em In celebration of the 125th}\\
{\em anniversary of the birth of}\\
{\em Ivan Matveevich Vinogradov}
\end{flushright}
\section{Introduction}

The familiar van der Corput $k$-th derivative estimate for exponential
sums (Titchmarsh \cite[Theorems 5.9, 5.11, \& 5.13]{Titch}, for
example), may be stated as follows. Let $k\ge 2$ be an integer, and
suppose that $f(x):[0,N]\to\R$ has continuous
derivatives of order up to $k$ on $(0,N)$.
Suppose further that $0<\lambda_k\le f^{(k)}(x)\le A\lambda_k$ on
$(0,N)$. Then
\beql{vdc}
\sum_{n\le N}e(f(n))\ll A^{2^{2-k}}N\lambda_k^{1/(2^k-2)}
+N^{1-2^{2-k}}\lambda_k^{-1/(2^k-2)},
\eeq
with an implied constant independent of $k$. One usually
chooses $k$ so that the first term dominates, and one often has
$A^{2^{2-k}}\ll 1$, so that the bound is merely
$O(N\lambda_k^{1/(2^k-2)})$. Clearly one can only get a non-trivial
bound when $\lk<1$.  A typical application is the series of
estimates
\[\zeta(\sigma+it)\ll t^{1/(2^k-2)}\log t,\;\;\;
\left(\sigma=1-\frac{k}{2^k-2},\; t\ge 2\right)\]
for $k=2,3,\ldots$. Again the implied constant is independent of $k$.

One can improve on the standard $k$-th derivative bound somewhat.  Thus
Robert and Sargos \cite{RS} show roughly that if $k=4$ then
\[\sum_{n\le N}e(f(n))\ll_{\ep} N^{\ep}(N\lambda_4^{1/13}+\lambda_4^{-7/13}),\]
for any $\ep>0$.  In the corresponding version of (\ref{vdc}) one
would have a term $N\lambda_4^{1/14}$ in place of
$N\lambda_4^{1/13}$. Similarly for $k=8$ and 9, Sargos \cite[Theorems
3 \& 4]{Sarg} gives bounds
\[\sum_{n\le N}e(f(n))\ll_{\ep}
N^{\ep}(N\lambda_8^{1/204}+\lambda_8^{-95/204}),\]
and
\[\sum_{n\le N}e(f(n))\ll_{\ep}
N^{\ep}(N\lambda_9^{7/2640}+\lambda_9^{-1001/2640}),\]
respectively.  Here the exponents $1/204$ and $7/2640$ should be
compared with the values $1/254$ and $1/510$ produced by (\ref{vdc}).
\bigskip

There are quite different approaches to exponential sums, using
estimates for the Vinogradov mean value integral
\beql{Jd}
J_{s,l}(P)=\int_0^1\ldots\int_0^1\left|\sum_{n\le P}e(\alpha_1
  n+\ldots +\alpha_l n^l)\right|^{2s}d\bal,
\eeq
see Vinogradov \cite{VB}, \cite{Vk}, and Korobov \cite{Kv}, amongst others.
The first of these methods is described by Titchmarsh \cite[Chapter 6]{Titch}
for example. The Vinogradov-Korobov machinery has been used by Ford
\cite[Theorem 2]{ford} to show that
\beql{f1}
\sum_{N<n\le 2N}n^{-it}\ll N^{1-1/134k^2}
\eeq
for $N^k\ge t\ge 2$.  (Ford's result is somewhat more precise, and
more general.)  One may think of this as corresponding very roughly to
a bound of the form (\ref{vdc}) with first term
$N\lambda_k^{1/134k^2}$.

A slightly refined version of
the original method of Vinogradov \cite{VB} coupled with 
new estimates for
the Vinogradov mean value integral,
leads to distinctly stronger
bounds.  For example, Wooley \cite[Theorem 1.2]{W3} gives
\[J_{s,l}(P)\ll_{\ep,l} P^{2s-l(l+1)/2+\ep}\;\;\;(s\ge l(l-1),\]
and Robert \cite[Theorem 10]{Rob} used this to show that if $k\ge 4$ then
\[\sum_{n\le N}e(f(n))\ll_{A,k,\ep} 
N^{1+\ep}(\lk^{1/2(k-1)(k-2)}+N^{-1/2(k-1)(k-2)})\]
for $N\ge\lk^{-(k-1)/(2k-3)}$.  This is a remarkable improvement on the
  classical $k$-th derivative estimate.  The exponent of $\lk$ is
  better than $1/(2^k-2)$ for all $k\ge 4$, and decreases
  quadratically rather than exponentially.

The purpose of this paper is to further refine
the original method of Vinogradov \cite{VB} and to input the very
recent optimal bounds for
the Vinogradov mean value integral, due to Wooley \cite{W3} (for $l=3$), and to
Bourgain, Demeter and Guth \cite{BDG} (for $l\ge 4$). These theorems
show that
\beql{bdg}
J_{s,l}(P)\ll_{\ep,l} P^{2s-l(l+1)/2+\ep}\;\;\;(s\ge\tfrac12 l(l+1),\;
l\ge 1),
\eeq
the cases $l=1$ and $l=2$ being elementary. The range for $s$ is
optimal, and it 
is this feature that represents the dramatic culmination of many
previous works over the past 80 years. Unfortunately neither result
gives an explicit dependence on $l$ and $s$, nor gives an explicit form
for the factor $P^{\ep}$. Results prior to the advent of Wooley's
efficient congruencing method had required $s$ to be
larger, but had given an explicit dependence on $l$.  Thus for
example, Ford
\cite[Theorem 3]{ford} implies in particular that
\[J_{s,l}(P)\ll l^{13l^3/4} P^{2s-l(l+1)/2+l^2/1000}\;\;\;
(s\ge\tfrac52 l^2,\; l\ge 129),\]
in which one has an additional term $l^2/1000$ in the exponent, and
more restrictive conditions on $l$ and $s$. An important application
of bounds for Weyl sums is to the zero-free region for $\zeta(s)$, as
described by Ford.  However for this it is crucial to have a suitable
dependence on the parameter $l$, so that the new result of Bourgain, 
Demeter and Guth is not applicable.

Our first result gives a new $k$-th derivative estimate
\begin{theorem}\label{esV}
Let $k\ge 3$ be an integer, and
suppose that $f(x):[0,N]\to\R$ has continuous
derivatives of order up to $k$ on $(0,N)$.  Suppose further that 
\[0<\lambda_k\le f^{(k)}(x)\le A\lambda_k,\;\;\; x\in (0,N). \]
Then
\[
\sum_{n\le N}e(f(n)) \ll_{A,k,\ep} 
N^{1+\ep}(\lk^{1/k(k-1)}+N^{-1/k(k-1)}+N^{-2/k(k-1)}\lk^{-2/k^2(k-1)}).
\]
\end{theorem}

If one thinks of $N\lk^{1/k(k-1)}$ as being the leading term here,
then one needs to compare the exponent $1/k(k-1)$ with the
corresponding exponent $1/(2^k-2)$ in (\ref{vdc}).  These agree for
$k=3$, but for larger values of $k$ the new exponent tends to zero far
more slowly than the old one.  It may perhaps be something of a
surprise that an analysis via Vinogradov's mean value integral
reproduces the same term $N\lambda_3^{1/6}$ as in the classical
third-derivative estimate.

We should emphasize that the stength of Theorem \ref{esV} comes almost
entirely from the new bound (\ref{bdg}). One could have injected
(\ref{bdg}) into the method of Robert \cite{Rob}, to produce an
estimate with the same terms $\lk^{1/k(k-1)}+N^{-1/k(k-1)}$ as in
Theorem \ref{esV}, but valid only for $N\ge \lk^{-(k-1)/(2k-3)}$.
Our result, incorporating a slightly better way of using the
Vinogradov mean value, gives the terms $\lk^{1/k(k-1)}+N^{-1/k(k-1)}$
in the substantially longer range $N\ge \lk^{-2/k}$.  However for our
application to Theorems \ref{EP}--\ref{rich} below, Robert's range
would have been very nearly sufficient.

The secondary terms in the bound given by Theorem \ref{esV} 
are somewhat awkward. The
classical estimate (\ref{vdc}) leads easily to an exponent pair, 
\[\left(\frac{1}{2^k-2}\,,\,\frac{2^k-k-1}{2^k-2}\right)\]
in which the term $N^{1-2^{2-k}}\lk^{-1/(2^k-2)}$ has no effect.
However the situation with Theorem \ref{esV} is more complicated. None 
the less we are able to produce a series of new exponent pairs. 

Before stating the
result we remind the reader of the necessary background. Let $s$ and
$c$ be positive constants, and let $\mathcal{F}(s,c)$ be the set of
quadruples $(N,I,f,y)$ where $y\ge N^s$ are positive real numbers, $I$
is a subinterval of $(N,2N]$, and $f$ is an infinitely differentiable
function on $I$, with 
\[\left|f^{(n+1)}(x)-\frac{d^n}{dx^n}(yx^{-s})\right|\le c
\left|\frac{d^n}{dx^n}(yx^{-s})\right|\]
for $x\in I$, for all $n\ge 0$.  We then say that $(p,q)$ is an
exponent pair, if $p$ and $q$ lie in the range
$0\le p\le\tfrac12\le q\le 1$, and for each $s$ there
is a corresponding $c=c(p,q,s)>0$ such that
\[\sum_{n\in I}e(f(n))\ll_{p,q,s} (yN^{-s})^p N^q,\]
uniformly for all quadruples $(N,I,f,y)\in\mathcal{F}(s,c)$.

We then have the following.
\begin{theorem}\label{EP}
For any integer $k\ge 3$ and any real $\ep>0$ there is an exponent 
pair given by
\beql{pqp}
p=\frac{2}{(k-1)^2(k+2)},
\eeq
and
\beql{pqq}
q=\frac{k^3+k^2-5k+2}{k(k-1)(k+2)}+\ep
=1-\frac{3k-2}{k(k-1)(k+2)}+\ep.
\eeq
\end{theorem}

In fact we are able to handle a much weaker condition on $f$. Let
$\mathcal{A}=(a_k)_3^\infty$ and $\mathcal{B}=(b_k)_3^\infty$ be
sequences of  positive real numbers, and let
$\mathcal{G}(\mathcal{A},\mathcal{B})$ be the set of
quadruples $(N,I,g,T)$ where $T\ge N$ are positive real numbers, $I$
is a subinterval of $(N,2N]$, and $g$ is an infinitely differentiable
function on $I$, with 
\[a_kTN^{-k}\le\left|g^{(k)}(x)\right|\le b_k TN^{-k}\]
for $x\in I$, and for all $k\ge 3$. We then have the following.
\begin{theorem}\label{EP0}
For any integer $k\ge 3$ and any real $\ep>0$, let $p$ and $q$ be given
by (\ref{pqp}) and (\ref{pqq}). Then
\[\sum_{n\in I}e(g(n))\ll_{k,\ep,\mathcal{A},\mathcal{B}} (TN^{-1})^p N^q,\]
uniformly for $(N,I,g,T)\in\mathcal{G}(\mathcal{A},\mathcal{B})$.
\end{theorem}
If $(N,I,g,y)\in\mathcal{F}(s,\tfrac14)$, then
$(N,I,g,yN^{1-s})\in\mathcal{G}(\mathcal{A},\mathcal{B})$ with
\[a_k=\frac{3\times 2^{1-2-k}}{4s(s+1)\ldots(s+k-2)},\;\;\;
b_k=\frac{5}{4s(s+1)\ldots(s+k-2)}.\]
The sequences $\mathcal{A}$ and $\mathcal{B}$ depend only on $s$, and
we immediately see that Theorem \ref{EP} follows from Theorem \ref{EP0}.

We next present a slightly weaker version of Theorem \ref{EP0}, which is
somewhat more immediately intelligible. It will be convenient to write
$T=N^{\tau}$. 
\begin{theorem}\label{EP1}
Let sequences $\mathcal{A}$ and $\mathcal{B}$, and a real number
$\ep>0$ be given, then
\[\sum_{n\in I}e(g(n))\ll_{\ep,\mathcal{A},\mathcal{B}} N^{1-49/(80\tau^2)+\ep},\]
uniformly for quadruples
$(N,I,g,T)\in\mathcal{G}(\mathcal{A},\mathcal{B})$ with $N\le T^{1/2}$.
\end{theorem}
The constant $49/80$ arises from the use of an exponent pair
\[(\tfrac{1}{20}\,,\,\tfrac{33}{40})=A^2BA^2B(0,1)\]
when $\tau=\tfrac72$. One could improve the constant slightly by employing a
better exponent pair. As will be clear from the proof, the constant
$\tfrac{49}{80}$ may be replaced by $1-\delta$ for any small
$\delta>0$, if we restrict to sufficiently large values
$\tau\ge\tau(\delta)$. 

As an example of Theorem \ref{EP1}, if $t\ge 2$ we find that
\beql{f2}
\sum_{n\in I}n^{-it}\ll_{\ep} N^{1-49/80\tau^2+\ep},
\eeq
for $\tau=(\log t)/(\log N)\ge 2$. This should be compared with
(\ref{f1}). Using (\ref{f2}) we produce the following result.
\begin{theorem}\label{rich}
Let $\kappa=\tfrac{8}{63}\sqrt{15}=0.4918\ldots$.  Then for any fixed
$\ep>0$ we have
\beql{rich1}
\zeta(\sigma+it)\ll_{\ep} t^{\kappa(1-\sigma)^{3/2}+\ep}
\eeq
uniformly for $t\ge 1$ and $\tfrac12\le\sigma\le 1$. Moreover we have
\beql{rich2}
\zeta(\sigma+it)\ll_{\ep} t^{\tfrac12 (1-\sigma)^{3/2}+\ep}
\eeq
uniformly for $t\ge 1$ and $0\le\sigma\le 1$.
\end{theorem}

One sees from the proof that $\kappa$ may be reduced to
$2/\sqrt{27}+\delta=0.3849\ldots$ for any small $\delta>0$, if we
restrict $\sigma$ to a suitably small range
$\sigma(\delta)\le\sigma\le 1$.
The corresponding result in the work of Ford \cite[Theorem 1]{ford}
states that
\[|\zeta(\sigma+it)|\le 76.2 t^{4.45(1-\sigma)^{3/2}}(\log t)^{2/3}\]
for $t\ge 3$ and $\tfrac12\le\sigma\le 1$.  Thus we have reduced the
constant $4.45$ to $0.4918\ldots$. Unfortunately our result yields no useful
information when $\sigma$ tends to 1, which is a critical situation in
many applications.  Moreover we do not have the explicit order constant
that Ford finds.

As Ford explains, there are a number of interesting corollaries, for which we
merely have to replace the constant $B=4.45$ by $B=0.492$ in the
arguments given in \cite[Pages 566 and 567]{ford}.  
We can feed our bound into the zero-density theorem of Montgomery
\cite[Theorem 12.3]{mont} (with $1-\alpha=4.93(1-\sigma)$ as used by Ford
\cite[Page 566]{ford}) to give the following.
\begin{corollary}
We have
\[N(\sigma,T)\ll_{\ep}T^{6.42(1-\sigma)^{3/2}+\ep}\]
for $\tfrac{9}{10}\le\sigma\le 1$.
\end{corollary}

For moments of the Riemann Zeta-function we have:
\begin{corollary}
For any positive integer $k$ one has
\[\int_{0}^T|\zeta(\sigma+it)|^{2k}dt\sim
T\sum_{1}^\infty d_k(n)^2n^{-2\sigma},\]
as $t\to\infty$, for any fixed $\sigma\ge 1-0.534 k^{-2/3}$.
\end{corollary}

For the generalized divisor problem we have:
\begin{corollary}
For any positive integer $k$ the error term $\Delta_k(x)$ in the
generalized divisor problem satisfies
\[\Delta(x)\ll_k x^{1-0.849 k^{-2/3}}.\]
\end{corollary}

In Section \ref{sec2} we will reduce the proof of Theorem \ref{esV} 
to a two-variable
counting problem involving fractional parts of the derivatives $f^{(j)}(n)$.
Section \ref{secN} shows how this counting problem is tackled, and
finally Section \ref{sPT} completes the proof of our theorems.

{\bf Acknowledgement.}  This work was supported by EPSRC grant number
EP/K021132X/1

\section{Initial Steps}\label{sec2}

Our goal in the first stage of the proof is to estimate the sum
\[\Sigma=\sum_{n\le N}e(f(n))\]
in terms of $J_{s,l}(P)$, together with, a counting function involving
the fractional parts of numbers of the form $f^{(j)}(n)/j!$.
\begin{lemma}\label{l1}
Let $k\ge 2$ be an integer, and
suppose that $f(x):[0,N]\to\R$ has continuous
derivatives of order up to $k$ on $(0,N)$.  Suppose further that 
\[0<\lambda_k\le f^{(k)}(x)\le A\lambda_k,\;\;\; x\in (0,N), \]
and that $A\lk\le \tfrac14$.  Then
\[\Sigma\ll H+k^2 N^{1-1/s}\mathcal{N}^{1/2s}
\left\{H^{-2s+k(k-1)/2}J_{s,k-1}(H)\right\}^{1/2s},\]
where $H=[(A\lk)^{-1/k}]$ and
\[\mathcal{N}=\#\left\{m,n\le N:
\left|\left|\frac{f^{(j)}(m)}{j!}-\frac{f^{(j)}(n)}{j!}\right|\right|
\le 2H^{-j}\mbox{ for } 1\le j\le k-1\right\}.\]
\end{lemma}

If $J_{s,k-1}(H)\ll_{\ep,k} H^{2s-k(k-1)/2+\ep}$ as in (\ref{bdg})
the estimate in the lemma reduces to
\beql{basic}
\Sigma\ll_{\ep,k} H+N^{1-1/s+\ep}\mathcal{N}^{1/2s}.
\eeq
Here we would want to choose $s$ to be as small as possible, and since
we are taking $l=k-1$ this means that we will have $s=k(k-1)/2$.

The lemma is clearly trivial if $H\ge N$, and we may therefore suppose
for the proof that $H\le N$.
For any positive integer $H\le N$ we will have
\[H\Sigma=\sum_{h\le H}\sum_{-h<n\le N-h}e(f(n+h))=
\sum_{h\le H}\sum_{1\le n\le N-H}e(f(n+h))+O(H^2),\]
so that
\beql{Sigma1}
\Sigma=H^{-1}\sum_{n\le N-H}\;\sum_{h\le H}e(f(n+h))+O(H).
\eeq

We proceed to approximate $f(n+h)$ by the polynomial
\[f_n(h):=f(n)+f'(n)h+\ldots+\frac{f^{(k-1)}(n)}{(k-1)!}h^{k-1}.\]
To do this we set $g_n(x)=f(n+x)-f_n(x)$ and use summation by parts to
obtain the bound
\[\sum_{h\le H}e(f(n+h))\ll |S_n(H)|+\int_0^H |S_n(x)g_n'(x)|dx,\]
where we have written
\[S_n(x)=\sum_{h\le x}e(f_n(h))\]
for convenience.

If $0\le x\le H$ we may use Taylor's Theorem with Lagrange's form of 
the remainder to show that
\[f'(n+x)=f_n'(x)+\frac{f^{(k)}(\xi)}{k!}x^{k-1}\]
for some $\xi\in(n,n+x)\subseteq(0,N)$. It follows that 
\[g_n'(x)\ll A\lk H^{k-1}\]
on $[0,H]$. With the choice $H=[(A\lk)^{-1/k}]$
we find that 
\[\sum_{h\le H}e(f(n+h))\ll |S_n(H)|+H^{-1}\int_0^H |S_n(x)|dx.\]
The bound (\ref{Sigma1}) now yields
\[\Sigma\ll H+H^{-1}\sum_{n\le N-H}|S_n(H)|+
H^{-2}\int_0^H \left\{\sum_{n\le N-H} |S_n(x)|\right\}dx.\]
It then follows that there is a positive integer $H_0\le H$ such that
\beql{b2}
\Sigma\ll H+H^{-1}\sum_{n\le N-H}|S_n(H_0)|.
\eeq

Now suppose that $\bal\in [0,1]^{k-1}$ and 
\beql{bal}
||f^{(j)}(n)/j!-\alpha_j||\le H^{-j}\mbox{ for } 1\le j\le k-1
\eeq
where 
\[||\theta||=\min_{n\in\Z}|\theta-n|\]
as usual. We proceed to replace $f_n(h)$ by
\[f(h;\bal)=\alpha_1 h+\ldots+\alpha_{k-1}h^{k-1}\]
as follows. Firstly we remove the constant term $f(n)$ from
$f_n(h)$. This has no effect on $|S_n(H_0)|$. Next, we replace each 
coefficient $f^{(j)}(n)/j!$ by $c_j$, say, with
$f^{(j)}(n)/j!-c_j\in\Z$, so that $|c_j-\alpha_j|\le H^{-j}$, and
denote the resulting polynomial by $f^*_n(h)$. If we write
\[S^*_n(H_0)=\sum_{h\le H_0}e(f^*_n(h))\]
then clearly $|S_n(H_0)|=|S^*_n(H_0)|$.  Moreover 
\[\frac{d}{dx}\left(f(x;\bal)-f^*_n(x)\right)\ll k^2\max_{j\le
  k-1}|c_j-\alpha_j|H^{j-1} \ll k^2H^{-1}.\]
It therefore follows on summing by parts that
\[S^*_n(H_0)\ll |S(H_0;\bal)|+k^2 H^{-1}\int_0^{H_0}|S(x;\bal)|dx,\]
where we have set
\[S(x;\bal)=\sum_{h\le x}e(f(h;\bal)).\]
We may therefore conclude that
\[S^*_n(H_0)\ll 2^{-k}H^{k(k-1)/2}\left\{\int_{\bal}|S(H_0;\bal)|
d\bal+k^2 H^{-1}\int_0^{H_0}\int_{\bal}|S(x;\bal)|
d\bal dx\right\},\]
where the integral over $\bal$ is for vectors in $[0,1]^{k-1}$
satisfying (\ref{bal}).

For each $\bal\in [0,1]^{k-1}$ we now define
\[\nu(\bal)=\#\{n\le N-H: ||f^{(j)}(n)/j!-\alpha_j||\le H^{-j}\mbox{
  for } 1\le j\le k-1\}.\]
We then find that
\beql{b1}
\sum_{n\le N-H}|S_n(H_0)|\ll 2^{-k}H^{k(k-1)/2}\left\{I(H_0)+
k^2 H^{-1}\int_0^{H_0}I(x)dx\right\},
\eeq
with
\[I(x)=\int_0^1\ldots\int_0^1|S(x;\bal)|\nu(\bal)
d\bal.\]
We easily see that
\[\int_0^1\ldots\int_0^1\nu(\bal)d\bal=
2^{k-1}H^{-k(k-1)/2}(N-H),\]
and that
\[\int_0^1\ldots\int_0^1\nu(\bal)^2d\bal\le
2^{k-1}H^{-k(k-1)/2}\mathcal{N},\]
where $\mathcal{N}$ is defined in Lemma \ref{l1}.
Moreover
\[\int_0^1\ldots\int_0^1|S(x;\bal)|^{2s}d\bal=J_{s,k-1}(x)\]
in the notation of (\ref{Jd}). Since $J_{s,k-1}(P)$ is non-decreasing
in $P$ this last integral may be bounded by $J_{s,k-1}(H)$.

Hence, by H\"{o}lder's inequality, for any positive integer $s$ we have
\[I(x)\ll
2^kH^{-k(k-1)/2}N^{1-1/s}\mathcal{N}^{1/2s}
\left\{H^{k(k-1)/2}J_{s,k-1}(H)\right\}^{1/2s}.\]
Thus (\ref{b1}) yields
\[\sum_{n\le N-H}|S_n(H_0)|\ll k^2 N^{1-1/s}\mathcal{N}^{1/2s}
\left\{H^{k(k-1)/2}J_{s,k-1}(H)\right\}^{1/2s}\]
and (\ref{b2}) gives us
\[\Sigma\ll H+k^2 N^{1-1/s}\mathcal{N}^{1/2s}
\left\{H^{-2s+k(k-1)/2}J_{s,k-1}(H)\right\}^{1/2s}\]
as required.

\section{The counting function $\mathcal{N}$}\label{secN}

Naturally our next task is to bound $\mathcal{N}$. The original
approach taken by Vinogradov, as described in Titchmarsh \cite[Chapter
6]{Titch}, merely used an $L^{\infty}$ bound for $\nu(\bal)$. One
discards all the information on $f^{(j)}(n)/j!$ for $j\le k-2$ and
uses only the case $j=k-1$. One then employs a standard procedure given by
the following trivial variant of \cite[Lemma 6.11]{Titch}, for
example.
\begin{lemma}\label{space}
Let $N$ be a positive integer, and suppose that $g(x):[0,N]\to\R$ has a
continuous derivative on $(0,N)$.  Suppose further that
\[0<\mu\le g'(x)\le A_0\mu,\;\;\; x\in(0,N).\]
Then
\[\#\{n\le N: ||g(n)||\le \theta\}\ll (1+A_0\mu N)(1+\mu^{-1}\theta).\]
\end{lemma}
We fix $m$ and take
\[g(x)=\frac{f^{(k-1)}(x)-f^{(k-1)}(m)}{(k-1)!}\]
and $\mu=\lk/(k-1)!$, $A_0=A$. This leads to a bound
\[\mathcal{N}\ll (k-1)!N(1+AN\lk)(1+H^{1-k}\lk^{-1}).\]
Under the assumption $A\lk\le \tfrac14$ in Lemma \ref{l1} we have
\[H^{1-k}\lk^{-1}\asymp (A\lk)^{1-1/k}\lk^{-1}=A(A\lk)^{-1/k}\ge A\ge 1,\]
whence our bound produces
\beql{addd}
\mathcal{N}\ll A^2(k-1)!N\lk^{-1/k}(1+N\lk).
\eeq
If one inserts this into (\ref{basic}) with $s=k(k-1)/2$ 
one gets an estimate
\begin{eqnarray*}
\Sigma&\ll_{\ep,k}&
(A\lk)^{-1/k}+N^{1-1/s+\ep}\{A^2N\lk^{-1/k}(1+N\lk)\}^{1/2s}\\
&\ll_{\ep,k}&
AN^{\ep}\{\lk^{-1/k}+N^{1-1/k(k-1)}\lk^{-1/k^2(k-1)}+N\lk^{1/k^2}\}.
\end{eqnarray*}
In fact the first term can be dropped, giving
\beql{vv}
\Sigma\ll_{\ep,k}AN^{\ep}\{N\lk^{1/k^2}+N^{1-1/k(k-1)}\lk^{-1/k^2(k-1)}\}.
\eeq
To see this we note that we have 
\[\Sigma\ll N^{1-1/k(k-1)}\lk^{-1/k^2(k-1)}\]
trivially
unless 
\[N^{1-1/k(k-1)}\lk^{-1/k^2(k-1)}\le N.\]
In this latter case however one sees that
\[\lk^{-1/k}\le N^{1-1/k(k-1)}\lk^{-1/k^2(k-1)}.\]

We may therefore regard (\ref{vv}) as being the result that
Vinogradov's method achieves, given the results of Wooley \cite{W3}
and Bourgain, Demeter and Guth \cite{BDG}.  It is already a remarkable
improvement on (\ref{vdc}), replacing the critical exponent
$1/(2^k-2)$ by $1/k^2$.  Thus, in appropriate circumstances, we get 
an improvement as soon as $k\ge 5$.  Our goal in this section is to
make the following small further sharpening in the estimation of
$\mathcal{N}$. 
\begin{lemma}\label{N2}
When $k\ge 3$ we have
\[\mathcal{N}\ll
\big((k-1)!A\big)^4(N+\lk N^2+\lk^{-2/k})\log N.\]
\end{lemma}
Apart from the term $\lk^{-2/k}$, which is insignificant in
applications, this represents an improvement of (\ref{addd}) by a
factor $\ll_{A,k}\lk^{1/k}$. 

On the one hand our proof will use the fact that $\mathcal{N}$ is a counting
function of two variables $m$ and $n$.  On the other we shall use information
about both $f^{(k-1)}$ and $f^{(k-2)}$.  The reader may find it
slightly surprising in the light of this that our bound depends on
$\lk$ only, and not on estimates for other derivatives $f^{(j)}$. The
introduction of $\mathcal{N}$, and our procedure for estimating it,
are the only really new aspects to this paper.
\bigskip

We begin our analysis by assuming that $k\ge 3$ and noting that
$\mathcal{N}$ is at most
\[\mathcal{N}_1=\#\left\{m,n\le N:
\left|\left|\frac{f^{(j)}(m)}{j!}-\frac{f^{(j)}(n)}{j!}\right|\right|
\le 2H^{-j}\mbox{ for } j=k-2, k-1\right\}.\]
We proceed to show that it suffices to consider pairs $m,n$ of integers that
are relatively close. It will be convenient to write $B=4H^{2-k}$ and
$C=4H^{1-k}$ and to set
\[g_1(x)=\frac{f^{(k-2)}(x)}{(k-2)!},\;\;\; g_2(x)=\frac{f^{(k-1)}(x)}{(k-1)!}.\]
We also define the doubly-periodic function
\[\phi(x,y)=\max\left(1-B^{-1}||x||,0\right)\max\left(1-C^{-1}||y||,0\right),\]
so that
\[\mathcal{N}_1\ll \sum_{m,n\le N}
\phi\big(g_1(m)-g_1(n)\,,\, g_2(m)-g_2(n)\big).\]
The function $\phi(x,y)$ has an absolutely convergent Fourier series
\[\phi(x,y)=\sum_{r,s\in\Z}c_{r,s}e(rx+sy)\]
with non-negative coefficients
\[c_{r,s}=BC\left(\frac{\sin(\pi rB)\sin(\pi sC)}{\pi^2
    rsBC}\right)^2.\]
Thus
\begin{eqnarray*}
\mathcal{N}_1&\ll& \sum_{r,s\in\Z}c_{r,s}\sum_{m,n\le N}
e\big(r(g_1(m)-g_1(n))+s(g_2(m)-g_2(n))\big)\\
&=& \sum_{r,s\in\Z}c_{r,s}\left|\sum_{n\le N}
e\big(rg_1(n)+sg_2(n)\big)\right|^2.
\end{eqnarray*}
Let $K$ be a positive integer parameter, to be chosen later.
We proceed to partition the range $(0,N]$ into $K$ intervals
$I_i=(a_i,b_i]$ for $i\le K$, having integer endpoints, and length
$b_i-a_i\le 1+N/K$.  An application of Cauchy's inequality then yields
\begin{eqnarray*}
\mathcal{N}_1 &\ll & K\sum_{i\le K}\sum_{r,s\in\Z}c_{r,s}
\left|\sum_{n\in I_i}e\big(rg_1(n)+sg_2(n)\big)\right|^2\\
&=&K\sum_{i\le K}\sum_{r,s\in\Z}c_{r,s}\sum_{m,n\in I_i}
e\big(r(g_1(m)-g_1(n))+s(g_2(m)-g_2(n))\big)\\
&=&K\sum_{i\le K}\sum_{m,n\in I_i}
\phi\big(g_1(m)-g_1(n)\,,\, g_2(m)-g_2(n)\big)\\
&\le & K\twosum{m,n\le N}{|m-n|\le 1+N/K}
\phi\big(g_1(m)-g_1(n)\,,\, g_2(m)-g_2(n)\big).
\end{eqnarray*}
We may therefore conclude that $\mathcal{N}_1\ll K\mathcal{N}_2$,
where $\mathcal{N}_2$ counts pairs of integers $m,n\le N$ with 
$|m-n|\le 1+N/K$ for which
\[\left|\left|\frac{f^{(j)}(m)}{j!}-\frac{f^{(j)}(n)}{j!}\right|\right|
\le 4H^{-j}\mbox{ for } j=k-2, k-1.\]

If $|m-n|\le 1+N/K$ we will have
\[\left|\frac{f^{(k-1)}(m)}{(k-1)!}-\frac{f^{(k-1)}(n)}{(k-1)!}\right|\le
\frac{|m-n|}{(k-1)!}\sup |f^{(k)}|\le A\lk(1+N/K),\]
by the mean-value theorem. We will choose
\[K=1+[4A\lk N],\]
so that
\[\left|\frac{f^{(k-1)}(m)}{(k-1)!}-\frac{f^{(k-1)}(n)}{(k-1)!}\right|\le
\frac{1}{2},\]
in view of our assumption that $A\lk\le\tfrac14$. Thus if
\[\left|\left|\frac{f^{(k-1)}(m)}{(k-1)!}-
\frac{f^{(k-1)}(n)}{(k-1)!}\right|\right|\le 4H^{1-k}\]
we must have
\[\left|\frac{f^{(k-1)}(m)}{(k-1)!}-\frac{f^{(k-1)}(n)}{(k-1)!}\right|
\le 4H^{1-k}.\]
However the mean-value theorem also tells us that
\[\left|\frac{f^{(k-1)}(m)}{(k-1)!}-\frac{f^{(k-1)}(n)}{(k-1)!}\right|\ge
\frac{|m-n|}{(k-1)!}\inf |f^{(k)}|\ge \lk\frac{|m-n|}{(k-1)!}.\]
We therefore conclude that
\[|m-n|\le \frac{4(k-1)!}{\lk H^{k-1}}\]
for any pair $m,n$ counted by $\mathcal{N}_2$.

There are $N$ pairs $m=n$ counted by $\mathcal{N}_2$. We consider
the remaining pairs with $m>n$, the alternative case producing the
same estimates by symmetry.  Then $m=n+d$ with $1\le d\le D$, where 
\[D=\min\left(N\,,\,\left[\frac{4(k-1)!}{\lk H^{k-1}}\right]\right).\]
For each available value of $d$ we estimate the number of
corresponding integers $n$ via Lemma \ref{space}, taking
\[g(x)=\frac{f^{(k-2)}(x+d)-f^{(k-2)}(x)}{(k-2)!}.\]
Then
\[g'(x)=\frac{f^{(k-1)}(x+d)-f^{(k-1)}(x)}{(k-2)!},\]
so that
\[d\frac{\lk}{(k-2)!}\le d\frac{\inf |f^{(k)}|}{(k-2)!}\le g'(x)\le
d\frac{\sup |f^{(k)}|}{(k-2)!}\le d\frac{A\lk}{(k-2)!},\]
by the mean-value theorem. We therefore apply the lemma with $\mu=
\lk d/(k-2)!$ and $A_0=A$. This shows that each $d\ge 1$ contributes
\begin{eqnarray*}
&\ll& (k-2)!(1+AN\lk d)(1+H^{2-k}\lk^{-1}d^{-1})\\
&\ll& (k-2)!(1+AN\lk D)(D+H^{2-k}\lk^{-1})d^{-1}\\
&\ll& \big((k-1)!\big)^3A(1+NH^{1-k})H^{2-k}\lk^{-1}d^{-1}\\
&\ll& \big((k-1)!A\big)^3(1+N\lk^{1-1/k})\lk^{-2/k}d^{-1}.
\end{eqnarray*}
Summing for $d\le D$ we therefore find that
\[\mathcal{N}_2\ll N+\big((k-1)!A\big)^3(1+N\lk^{1-1/k})\lk^{-2/k}\log D.\]
Since $k\ge 3$, $\lk\le 1$ and $D\le N$ this simplifies to give
\[\mathcal{N}_2\ll \big((k-1)!A\big)^3(N+\lk^{-2/k})\log N,\]
whence
\[\mathcal{N}\le\mathcal{N}_1\ll K\mathcal{N}_2\ll (1+A\lk N)
\big((k-1)!A\big)^3(N+\lk^{-2/k})\log N\]
\[\ll \big((k-1)!A\big)^4(N+\lk N^2+\lk^{-2/k}+N\lk^{1-2/k})\log N.\]
Since $k\ge 3$ and $\lk\le 1$ we have $N\lk^{1-2/k}\le N$, and 
Lemma \ref{N2} follows.

\section{Proof of the Theorems}\label{sPT}

If we insert Lemma \ref{N2} into Lemma \ref{l1}, and use the bound
(\ref{bdg}) with the choices $l=k-1, s=k(k-1)/2$, we see that
\[\Sigma\ll_{A,k,\ep}
N^{\ep}(\lk^{-1/k}+N^{1-1/k(k-1)}+N\lk^{1/k(k-1)}
+N^{1-2/k(k-1)}\lk^{-2/k^2(k-1)}).\]
The term $\lk^{-1/k}$ may be omitted, since the resulting bound
\[\Sigma\ll_{A,k,\ep}N^{\ep}(N^{1-1/k(k-1)}+N\lk^{1/k(k-1)}
+N^{1-2/k(k-1)}\lk^{-2/k^2(k-1)})\]
holds trivially when $N\le N^{1-2/k(k-1)}\lk^{-2/k^2(k-1)}$, while
\[\lk^{-1/k}\le N^{1-2/k(k-1)}\lk^{-2/k^2(k-1)}\]
when $N\ge N^{1-2/k(k-1)}\lk^{-2/k^2(k-1)}$.
This suffices for Theorem \ref{esV}.

We turn next to Theorem \ref{EP0}.  Suppose that 
$(N,I,g,T)\in\mathcal{G}(\mathcal{A},\mathcal{B})$, and let $I$ have end points
$N_0$ and $N_0+N_1$, so that $N_1\le N$. We apply Theorem \ref{esV}
to the function $f(x)=g(N_0+x)$, taking $\lambda_k=a_kTN^{-k}$ and
$A=b_k/a_k$. (Since $f^{(k)}$ is differentiable it is continuous, and
hence it cannot change sign if $|f^{(k)}(x)|\ge a_kTN^{-k}>0$. Taking
complex conjugates of our sum if necessary we may therefore assume
that $f^{(k)}(x)$ is positive on $I$.) It follows that if $k\ge 3$ then
\begin{eqnarray*}
\lefteqn{\sum_{n\in I}e(g(n))}\\
&\ll_{\ep,k,\mathcal{A},\mathcal{B}}&
N^{1+\ep}(\lambda_k^{1/k(k-1)}+N^{-1/k(k-1)}+N^{-2/k(k-1)}\lambda_k^{-2/k^2(k-1)})\\
&\ll_{\ep,k,\mathcal{A},\mathcal{B}}&
N^{1+\ep}(N^{-1/(k-1)}T^{1/k(k-1)}+N^{-1/k(k-1)}+T^{-2/k^2(k-1)}).
\end{eqnarray*}

We use the above bound for 
\[\frac{(k-1)^2+1}{k}\le\tau<\frac{k^2+1}{k+1}\]
where we define $\tau$ by  $T=N^{\tau}$.  For this range of $\tau$ we
find that
\begin{eqnarray}\label{mid}
\lefteqn{\max\left(\frac{\tau-k}{k(k-1)}\,,\,\frac{-1}{k(k-1)}\,,\,
\frac{-2\tau}{k^2(k-1)}\right)}\nonumber\\
&=&\left\{\begin{array}{cc}
-1/k(k-1), & ((k-1)^2+1)/k \le \tau\le k-1,\\
(\tau-k)/k(k-1), & k-1\le\tau<(k^2+1)/(k+1),\end{array}\right.\\
&\le & A_k\tau+B_k,\nonumber
\end{eqnarray}
where the coefficients $A_k$ and $B_k$ are chosen so that
\[A_k\frac{(k-1)^2+1}{k}+B_k=\frac{-1}{k(k-1)}\]
and
\[A_k\frac{k^2+1}{k+1}+B_k=\frac{(k^2+1)/(k+1)-k}{k(k-1)}=\frac{-1}{k(k+1)}.\]
One then calculates that 
\[A_k=\frac{2}{(k-1)^2(k+2)}\;\;\;\mbox{and}\;\;\;
B_k=-\frac{3k^2-3k+2}{k(k-1)^2(k+2)}.\]
If we now define $\phi(\tau):[2,\infty)\to\R$ by taking
$\phi(\tau)=A_k\tau+B_k$ on 
\[\left[\frac{(k-1)^2+1}{k}\,,\,\frac{k^2+1}{k+1}\right)\]
for each integer $k\ge 3$, we conclude that
\beql{nd}
\sum_{n\in I}e(g(n))\ll_{\ep,\tau_0,\mathcal{A},\mathcal{B}}
N^{1+\phi(\tau)+\ep},
\eeq
uniformly for $2\le\tau\le\tau_0$. The function $\phi$ is continuous,
and since the coefficients $A_k$ are monotonic decreasing $\phi$ is
also convex. It follows that $\phi(\tau)\le A_k\tau+B_k$ for any
$\tau\in[2,\infty)$ and any $k\ge 3$. Thus
\[\sum_{n\in I}e(g(n))\ll_{\ep,\tau_0,\mathcal{A},\mathcal{B}}
N^{1+B_k+\ep}T^{A_k}=(TN^{-1})^p N^q,\]
with $p,q$ given by (\ref{pqp}) and (\ref{pqq}). 
As before, this is uniform in any
finite range $2\le\tau\le\tau_0$.  However if we set $\tau_0=1+(1-q)/p$
then $\tau_0$ will depend on $\ep$ and $k$ alone.  Moreover, if
$\tau\ge\tau_0$ then we trivially have
\[\sum_{n\in I}e(g(n))\ll N\le (TN^{-1})^pN^q.\]
Finally, if $\tau\le 2$ we use the well known exponent pair
$(\tfrac16,\tfrac23)$ to show that
\[\sum_{n\in I}e(g(n))\ll T^{1/6}N^{1/2}.\]
When $k\ge 3$ one easily verifies that $q\ge p+1/2$ and $p+q\ge 5/6$
for the values (\ref{pqp}) and (\ref{pqq}), whence $T^{1/6}N^{1/2}\le
T^pN^{q-p}$ for $N\ge T^{1/2}$. It then follows that
\[\sum_{n\in I}e(g(n))\ll T^{1/6}N^{1/2} \le T^pN^{q-p}\]
for the remaining range $1\le\tau\le 2$.  This completes the proof of
Theorem \ref{EP0}.
\bigskip

We move now to the proof of Theorem \ref{EP1}. Let
$\tau_0=\sqrt{49/80\ep^2}$.  Then if $\tau\ge\tau_0$ we will trivially
have
\[\sum_{n\in I}e(g(n))\ll N\le N^{1-49/80\tau^2+\ep}.\]
When $\tau\le\tau_0$ we begin by handling
the range $\tfrac{13}{3}\le\tau\le\tau_0$, for which we claim that
$\phi(\tau)\le -49/80\tau^2$. This will clearly suffice, in view of the
estimate (\ref{nd}). Since $\phi(\tau)$ is piecewise linear, while the
function $-49/80\tau^2$ is convex, it suffices to verify that
$\phi(\tau)\le -49/80\tau^2$ at each of the points
$\tau=(k^2+1)/(k+1)$, for $k\ge 5$. This condition is equivalent to
\[\frac{(k^2+1)^2}{k(k+1)^3}\ge\frac{49}{80}.\]
However the fraction on the right is increasing for $k\ge 5$, and
takes the value $169/270>49/80$ at $k=5$.

When $\tfrac72\le\tau\le\tfrac{13}{3}$ we will use the bounds
\[\sum_{n\in I}e(g(n))\ll_{\ep} N^{1-1/20+\ep},\;\;\;(\tfrac72\le\tau\le
4)\]
and
\[\sum_{n\in I}e(g(n))\ll_{\ep} N^{1-(5-\tau)/20+\ep},\;\;\;(4\le\tau\le
\tfrac{13}{3})\]
which come from the case $k=5$ of (\ref{mid}). Note that the first of
these is valid in the longer range $\tfrac{17}{5}\le\tau\le 4$, but we
shall only use it when $\tfrac72\le\tau\le 4$. We therefore need to
verify that $-1/20\le-49/80\tau^2$ for $\tfrac72\le\tau\le 4$ and that
$-(5-\tau)/20\le-49/80\tau^2$ for $4\le\tau\le \tfrac{13}{3}$. This
is routine, but we observe that we have equality at $\tau=\tfrac72$.

We next consider the case in which $\tfrac{59}{22}\le\tau\le\tfrac72$, 
for which we use the bound
\[\sum_{n\in I}e(g(n))\ll (T/N)^{1/20}N^{33/40}=N^{1+(2\tau-9)/40}\]
corresponding to the exponent pair
$(\tfrac{1}{20},\tfrac{33}{40})$. (This pair is $A^2BA^2B(0,1)$ in the
usual notation, see Titchmarsh \cite[\S 5.20]{Titch}, for example.)
Again, it is routine to check that
\[\frac{2\tau-9}{40}\le-\frac{49}{80\tau^2},\;\;\;(\frac{59}{22}\le
\tau\le\frac{7}{2}).\]

Finally we examine the range $2\le\tau\le\tfrac{59}{22}$, and here we
use the bound
\[\sum_{n\in I}e(g(n))\ll (T/N)^{1/9}N^{13/18}=N^{1+(2\tau-7)/18}\]
corresponding to the exponent pair
$(\tfrac{1}{9},\tfrac{13}{18})$. (This pair is $ABA^2B(0,1)$ in the
usual notation, see Titchmarsh \cite[\S 5.20]{Titch}, for example.)
Another routine check shows that
\[\frac{2\tau-7}{18}\le-\frac{49}{80\tau^2},\;\;\;(2\le
\tau\le\frac{59}{22}),\]
thereby completing the proof of Theorem \ref{EP1}.
\bigskip

We turn now to Theorem \ref{rich}. If $\tau\ge 2$ we may use 
(\ref{f2}) along with a
partial summation to obtain
\[\sum_{n\in J}n^{-\sigma-it}\ll_{\ep} N^{1-49/80\tau^2-\sigma+\ep}\le
t^{(1-\sigma)\tau^{-1}-\tfrac{49}{80}\tau^{-3}+\ep/2}\]
for any $\sigma\in[\tfrac12,1]$, and for any interval $J\subseteq(N,2N]$.
As a function of
$\tau\in(0,\infty)$ the exponent of $t$ is maximal at
\[\tau=\sqrt{\frac{147}{80(1-\sigma)}},\]
whence
\[\sum_{n\in J}n^{-\sigma-it}\ll_{\ep}t^{\kappa(1-\sigma)^{3/2}+\ep/2}.\]
Using a dyadic subdivision of $(0,N]$ we therefore have
\[\sum_{n\le N}n^{-\sigma-it}\ll_{\ep}t^{\kappa(1-\sigma)^{3/2}+3\ep/4}\]
for any $N\le t^{1/2}$. A further summation by parts then shows that
\[\sum_{n\le M}n^{-1+\sigma-it}
\ll_{\ep}M^{2\sigma-1}t^{\kappa(1-\sigma)^{3/2}+\ep}
\ll_{\ep}t^{\sigma-\tfrac12+\kappa(1-\sigma)^{3/2}+\ep}\]
for any $M\le t^{1/2}$. The required bound (\ref{rich1}) then follows
from the approximate functional equation for $\zeta(s)$.

The bound (\ref{rich2}) follows from (\ref{rich1}) when
$\tfrac12\le\sigma\le 1$, since $\kappa<\tfrac12$. For the remaining
range we use the functional equation, which shows that
\[\zeta(\sigma+it)\ll t^{\tfrac12-\sigma}|\zeta(1-\sigma+it)|\ll_{\ep}
t^{\tfrac12-\sigma+\tfrac12 \sigma^{3/2}+\ep}.\]
However one can readily verify that
\[\frac{1}{2}-\sigma+\frac{\sigma^{3/2}}{2}\le\frac{(1-\sigma)^{3/2}}{2}\]
for $0\le\sigma\le\tfrac12$, which completes the proof of Theorem
\ref{rich}.

\bigskip
\bigskip

Mathematical Institute,

Radcliffe Observatory Quarter,

Woodstock Road,

Oxford

OX2 6GG

UK
\bigskip

{\tt rhb@maths.ox.ac.uk}


\begin{thebibliography}{99}
\bibitem{BDG}J. Bourgain, C. Demeter and L. Guth, Proof of the main
  conjecture in Vinogradov's mean value theorem for degrees higher
  than three,  {\tt arXiv:1512.01565}.
\bibitem{ford} K. Ford, Vinogradov's integral and bounds for 
the Riemann zeta function, {\em Proc. London Math. Soc. (3)}, 
85 (2002), no. 3, 565--633. 
\bibitem{Kv} N.M. Korobov, Estimates of trigonometric sums and 
their applications, {\em Uspehi Mat. Nauk}, 13 (1958) no. 4 (82), 185--192.
\bibitem{mont}H.L. Montgomery, {\em Topics in multiplicative number
theory,} Lecture Notes in Mathematics 227, (Springer, Berlin, 1971).
\bibitem{Rob}O. Robert, On van der Corput’s $k$-th derivative test
  for exponential sums, {\em Indag. Math. (N.S.)}, to appear.
\bibitem{RS} O. Robert and P. Sargos, A fourth derivative test for 
exponential sums, {\em Compositio Math.}, 130 (2002), no. 3, 275--292.
\bibitem{Sarg} P. Sargos, An analog of van der Corput's $A^4$-process 
for exponential sums, {\em Acta Arith.}, 110 (2003), no. 3, 219--231.
\bibitem{Titch} E.C. Titchmarsh, {\em The theory of the Riemann 
zeta-function}, Second edition, (Clarendon Press, Oxford University
Press, New York, 1986).
\bibitem{VB} I.M. Vinogradov, New Estimates for Weyl Sums, {\em
    Doklady Nauk SSSR}, 8 (1935), no. 5, 195--198.
\bibitem{Vk} I.M. Vinogradov, A new estimate of the function $\zeta(1+it)$,
{\em Izv. Akad. Nauk SSSR. Ser. Mat.}, 22 (1958), 161--164. 
\bibitem{W3} T.D. Wooley, The cubic case of the main conjecture in 
Vinogradov's mean value theorem, {\tt arXiv:1401.3150}.

\end{thebibliography}
\end{document}